\documentclass[a4paper,12pt]{article}
\usepackage[english]{babel}
\usepackage{hyperref}
\usepackage{indentfirst}
\usepackage{amsfonts,amssymb,amsthm,amsmath}
 \theoremstyle{definition}
\newtheorem{De}{Definition}
\theoremstyle{theorem}
\newtheorem{Th}{Theorem}
\theoremstyle{lemma}

\theoremstyle{plain}
\newtheorem{Pl}{Proposition}
\theoremstyle{remark}
\newtheorem{Remark}{Remark}
\theoremstyle{plain}
\newtheorem{coll}{Corollary}
\newcounter{tmp}
\newcommand{\ir}{[t_0,\vartheta_0]\times \mathbb{R}^n}
\newcommand{\Dh}[2]{{\rm D}_H{#1}(#2)}
\begin{document}
\title{Characterization of Feedback Nash Equilibrium for Differential Games\footnote{The work is supported by RFBR (grant No 09-01-00436-a), Grant of President of Russian Federation (project MK-7320.2010.1), RAS Presidium Program of Fundamental Researches ``Mathematical Theory of Control''.}}
\author{Yurii Averboukh}
\date{\fontsize{10}{12}\selectfont Institute of Mathematics and Mechanics UrB RAS\\
S. Kovalevskaya street 16\\
620219, GSP-384, Ekaterinburg
Russia\\
ayv@imm.uran.ru, averboukh@gmail.com}
\maketitle
\abstract{We investigate the set of Nash equilibrium payoffs for two person differential games. The main result of the paper is the characterization of the set of Nash equilibrium payoffs in the terms of nonsmooth analysis. Also we obtain the sufficient conditions for a pair of continuous function to provide a Nash equilibrium.  This result generalizes the method of system of Hamilton-Jacobi equations. }
\section{Introduction}
In this paper we characterize Nash equilibrium payoffs for two person differential games. We consider nonzero-sum differential games in the framework of positional strategies first suggested by N.N. Krasovskii for zero-sum differential games \cite{NN_PDG_en}. The existence of Nash equilibrium was established in works of A.F. Kononenko \cite{Kononenko} and A.F. Kleimenov \cite{Cleimenov}. The proof is based on punishment strategies technique.
This technique permits to characterize the set of Nash equilibrium payoffs \cite{Cleimenov}, \cite{Chistyakov_Nash}.

The main result of this paper is the characterization of the set of Nash equilibrium payoffs in the terms of nonsmooth analysis. Also we obtain the sufficient conditions for a pair of continuous function to provide a Nash equilibrium.  This result generalizes the method of the systems of Hamilton-Jacobi equations.

\section{Preliminaries}
We consider the following doubly controlled system
\begin{equation}\label{sys}
 \dot{x}=f(t,x,u,v), \ \ t\in [t_0,\vartheta_0], \ \ x\in\mathbb{R}^n, \ \ u\in P,\ \  v\in Q.
\end{equation} Here $u$ and $v$ are controls of the player I and the player II respectively. Payoffs are terminal. The player I wants to maximize $\sigma_1(x(\vartheta_0))$, the player II wants to maximize $\sigma_2(x(\vartheta_0))$. We assume that the sets $P$ and $Q$ are compacts, the function $f$, $\sigma_1$ and $\sigma_2$ are continuous, moreover $f$ is Lipschitz continuous with respect to the phase variable, and satisfies the sublinear growth condition with respect to $x$.

The use the control design suggested in \cite{Cleimenov}. This control design follows N.N. Krasovskii positional formalization. Feedback strategy of the Player I is a pair of function $U=(u(t,x,\varepsilon),\beta_1(\varepsilon))$. Here $u(t,x,\varepsilon)$ is a function of position $(t,x)\in \ir$ and precision parameter $\varepsilon$, $\beta_1(\varepsilon)$ is a continuous function of precision parameter. We suppose that $\beta_1(\varepsilon)\rightarrow 0,$ $\varepsilon\rightarrow 0$. Analogously, the feedback strategy of the Player II is a pair $V=(v(t,x,\varepsilon),\beta_2(\varepsilon))$.

Let a position $(t_*,x_*)$ be chosen. Step-by-step motion is defined in the following way. We suppose that the players choose precision parameters  $\varepsilon_1$ and $\varepsilon_2$ respectively. Let the Player I choose the partition of the interval  $[t_*,\vartheta_0]$ $\Delta_1=\{\tau_j\}_{j=0}^{r}$ of the fineness less than $\varepsilon_1. $ Suppose that the Player II chooses the partition  $\Delta_2=\{\xi_k\}_{k=1}^\nu$ of the fineness less than  $\varepsilon_2$. The solution $x[\cdot]$ of equation (\ref{sys}) with initial date $x[t_*]=x_*$ such that the control of the Player I is equal to $u(\tau_j,x[\tau_j],\varepsilon_1)$ on $[\tau_j,\tau_{j+1})$, and the control of the Player II is  equal to  $v(\xi_k,x[\xi_k],\varepsilon_2)$ on $[\xi_k,\xi_{k+1})$ is called a step-by-step motion. Denote it by $x[\cdot,t_*,x_*;U,\varepsilon_1,\Delta_1;V,\varepsilon_2;\Delta_2]$. The set of all step-by-step motions from the position $(t_*,x_*)$ under strategies $U$ and $V$ and precision parameters $\varepsilon_1$ and $\varepsilon_2$ is denoted by $X(t_*,x_*;U,\varepsilon_1;V,\varepsilon_2)$. The step-by-step motions is called consistent if $\varepsilon_1=\varepsilon_2$.

A limit of step-by-motions $x[\cdot,t^k,x^k;U,\varepsilon_1^k,\Delta^k_1;V,\varepsilon^k_2,\Delta^k_2]$ is called constructive motion if $t^k\rightarrow t_*$, $x^k\rightarrow x_*$, $\varepsilon^k_1\rightarrow 0$, $\varepsilon^k_2\rightarrow 0$, $k\rightarrow \infty$. Denote by $X(t_*,x_*;U,V)$ the set of constructive motions. By Arzela-Ascoli theorem the set of constructive motions is nonempty. If the limit is taken only by consistent step-by-step motions the limit is called consistent constructive motions. Denote the set of consistent constructive motions by $X^c(t_*,x_*;U,V)$. This set is nonempty also.

The following definition of Nash equilibrium is used.
\begin{De}Let $(t_*,x_*)\in \ir$. The pair of strategies $U^N$ and $V^N$ is said to be Nash equilibrium solution at the position $(t_*,x_*)$, if for all strategies $U$ and $V$ the following inequalities hold:
$$
    \max\{\sigma_1(x[\vartheta_0]):x[\cdot]\in X(t_*,x_*,U,V^N)\}\leq  \min\{\sigma_1(x^c[\vartheta_0]):x^c[\cdot]\in X^c(t_*,x_*,U^N,V^N)\}.
$$
$$
    \max\{\sigma_2(x[\vartheta_0]):x[\cdot]\in X(t_*,x_*,U^N,V)\}\leq  \min\{\sigma_2(x^c[\vartheta_0]):x^c[\cdot]\in X^c(t_*,x_*,U^N,V^N)\}.
$$

\end{De}

The pair of payoff $(J_1,J_2)$ determined by a Nash solution is called a Nash equilibrium payoff of the game. In the typical case there are many Nash equilibriums with different payoffs. The set of all  Nash equilibrium payoffs is called a Nash value of the game and is denoted by by $\mathcal{N}(t_*,x_*)$. One can consider multivalued map taking   $(t_*,x_*)$ to the $\mathcal{N}(t_*,x_*)$.

The set $\mathcal{N}(t_*,x_*)$ is nonempty under the Isaacs condition \cite{Cleimenov}, \cite{Kononenko}. The proof is based on the punishment strategy technique. If the Isaacs condition is not fulfilled the Nash equilibrium solution exists in the class of mixed strategies or in the class of pair counterstrategy/strategy \cite{Cleimenov}.

Below we suppose the Isaacs condition holds:
for all $t\in [t_0,\vartheta_0],x,s\in\mathbb{R}^n$ $$\min_{u\in P}\max_{v\in Q}\langle s, f(t,x,u,v)\rangle=\max_{v\in Q}\min_{u\in P}\langle s, f(t,x,u,v)\rangle. $$

\begin{Remark} If the Isaacs condition doesn't hold, one can consider the solution in the class of mixed strategies. For this purpose we consider doubly controlled system
 \begin{equation}\label{mixed_system}
 \dot{x}=\int_{P}\int_{Q}f(t,x,u,v)\nu(dv)\mu(du), \ \ t\in [t_0,\vartheta], \ \ x\in\mathbb{R}^n, \ \ \mu\in{\rm rpm}(P), \ \ \nu\in{\rm rpm}(Q).
\end{equation}
Here $\mu$ is a generalized control of the Player I, $\nu$ is a generalized control of the Player II, ${\rm rpm}(P), {\rm rpm}(Q)$ are sets of regular probabilistic measures on $P$ and $Q$ respectively. We endow the sets  ${\rm rpm}(P)$ and ${\rm rpm}(Q)$ with $*$-weak topology. Obtained topology spaces are compacts. It is easy to show that the Isaacs condition is fulfilled for system (\ref{mixed_system}). Further we will not mention  the change from system (\ref{sys}) to system (\ref{mixed_system}).
\end{Remark}

Consider the zero-sum differential game $\Gamma_1$ with dynamic determined by (\ref{sys}) and the payoff determined by $\sigma_1(x(\vartheta_0))$. We assume that the Player I wants to maximize  $\sigma_1(x(\vartheta_0))$ the interest of the Player II is opposite. There exists the value of the game $\Gamma_1$. Denote it by $\omega_1$. Analogously consider the zero-sum differential game with the dynamics (\ref{sys}) and the payoff $\sigma_2$. We assume that the Player II wants to maximize $\sigma_2(x(\vartheta_0))$ while the Player I want to minimize it. Denote the value of this game by $\omega_2$.

\section{Main result}\label{sec_main}
Consider the differential inclusion \begin{equation}\label{F_diff_incl_def}
\dot{x}\in \mathcal{F}(t,x)\triangleq {\rm co}\{f(t,x,u,v):u\in P,v\in Q\}. \end{equation} By ${\rm Sol}(t_*,x_*)$ denote the set of solution of (\ref{F_diff_incl_def}) with initial data  $x(t_*)=x_*$.
\begin{Pl}\label{th_nash_inclusion}
Let the multivalued map $\mathcal{T}:\ir\rightarrow \mathcal{P}(\mathbb{R}^2)$ satisfy the following conditions:
\begin{list}{\rm (N\arabic{tmp})}{\usecounter{tmp}}
 \item $\mathcal{T}(t,x)\subset [\omega_1(t,x),\infty)\times[\omega_2(t,x),\infty)$ for all $(t,x)\in \ir$;
 \item $\mathcal{T}(\vartheta_0,x)=\{(\sigma_1(x),\sigma_2(x))\}$ for all $x\in\mathbb{R}^n$;
 \item for all $(t_*,x_*)\in\ir$, $(J_1,J_2)\in\mathcal{T}(t_*,x_*)$ there exists a motion ${y(\cdot)\in {\rm Sol}(t_*,x_*)}$ such that
 $$(J_1,J_2)\in \mathcal{T}(t,y(t)), \ \ t\in [t_*,\vartheta_0].$$
\end{list} Then $\mathcal{T}(t,x)\subset\mathcal{N}(t,x)$ for all $(t,x)\in \ir$.
\end{Pl}
This proposition follows from \cite[Theorem 1.4]{Cleimenov}.

Further we limit our attention to closed multivalues maps. The  map $\mathcal{T}:\ir\rightarrow \mathcal{P}(\mathbb{R}^2)$ is called closed if its graph is closed $\mathcal{T}$, i.e. ${\rm Cl}[\mathcal{T}]=\mathcal{T}$. Here ${\rm Cl}$ denotes the closure of graph: \begin{multline*}[{\rm Cl} \mathcal{T}](t,x)\triangleq \{(J_1,J_2):\exists\{(t^k,x^k)\}_{k=1}^\infty\subset\ir\ \ \exists\{(z^k_1,z^k_2)\}\subset\mathbb{R}^2:\\ (z_1^k,z^k_2)\in \mathcal{T}(t^k,x^k), \ \  (t^k,x^k)\rightarrow (t,x), \ \ (z^k_1,z_2^k)\rightarrow (J_1,J_2), \ \ k\rightarrow\infty\}. \end{multline*}

Let $I$ be a indexing set. Let multivalued maps  $\mathcal{T}^\alpha:\ir\rightarrow \mathcal{P}(\mathbb{R}^2)$, $\alpha\in I$, satisfy  conditions (N1)--(N3). Define the map $\mathcal{T}^*:\ir\rightarrow \mathcal{P}(\mathbb{R}^2)$ by the rule $\mathcal{T}^*={\rm Cl}\mathcal{T}$, where
$$\mathcal{T}(t,x)\triangleq \Bigl[\bigcup_{\alpha\in I}\mathcal{T}^\alpha(t,x)\Bigr]. $$ The multivalued map $\mathcal{T}^*$ is closed has compact images and satisfies conditions (N1)--(N3).
By $\mathcal{T}^+$ denote the closure of pointwise union of all upper semicontinuous  multivalued map from  $\ir$ to  $\mathbb{R}^2$ satisfying conditions  (N1)--(N3).
It follows  \cite{Cleimenov} that $\mathcal{T}^+(t,x)=\mathcal{N}(t,x)$ for all
 $(t,x)\in \ir$.

Further we formulate  condition (N3) in the terms of viability theory and obtain the infinitesimal form of this condition.
\begin{Th}\label{th_weak_inv} Let the map $\mathcal{T}:\ir\rightarrow {\mathcal{P}}(\mathbb{R}^2)$ be closed.
Then  condition {\rm (N3)} is equivalent to the following one: for all $(t_*,x_*)\in \ir$, $(J_1,J_2)\in\mathcal{T}(t_*,x_*)$ there exist $\theta>t_*$ and $y(\cdot)\in {\rm Sol}(t_*,x_*)$ such that
$$
(J_1,J_2)\in\mathcal{T}(t,y(t)), \ \ t\in [t_*,\theta].
$$
\end{Th}
Theorem \ref{th_weak_inv} is proved in section \ref{sec_weak}.

In order to obtain the infinitesimal form of condition (N3) we define a derivative of a multivalued map. By ${\rm dist}$ denote the following planar distance between the point $(J_1,J_2)\in\mathbb{R}^2$ and the set $A\subset\mathbb{R}^2$:
$${\rm dist}[(J_1,J_2),A]\triangleq  \inf\{|\zeta_1-J_1|+|\zeta_2-J_2|:(\zeta_1,\zeta_2)\in A \}.$$
Define the directional derivative of the multivalued map by the rule $$
\Dh{\mathcal{T}}{t,x; (J_1,J_2),w}\triangleq\liminf_{\delta\downarrow 0,w'\rightarrow w}\frac{{\rm dist}[(J_1,J_2),\mathcal{T}(t+\delta, x+\delta w')]}{\delta}. $$
\begin{Th}\label{th_infinitizemal} Let $\mathcal{T}:\ir\rightarrow \mathcal{P}(\mathbb{R}^2)$ be closed.
Then condition {\rm (N3)} at the position $(t_*,x_*)\in \ir$ is equivalent to the following one:    \begin{equation}\label{dh_zero}
\sup_{(J_1,J_2)\in\mathcal{T}(t_*,x_*)}\inf_{w\in\mathcal{F}(t_*,x_*)}\Dh{\mathcal{T}}{t_*,x_*;(J_1,J_2),w}=0.
\end{equation}
\end{Th}
Theorem \ref{th_infinitizemal} is proved in section \ref{sec_weak}.

Introduce the set $$\widehat{\partial}\mathcal{T}(t_*,x_*;(J_1;J_2))\triangleq \left\{w:\Dh{\mathcal{T}}{t_*,x_*;(J_1,J_2),w}=0\right\}.$$

\begin{Remark}\label{rm_diff}
Condition (\ref{dh_zero}) can be formulated in the following way:
$$\widehat{\partial}\mathcal{T}(t_*,x_*;(J_1;J_2))\cap \mathcal{F}(t_*,x_*)=0,\ \ \forall (J_1,J_2)\in\mathcal{T}(t_*,x_*). $$\end{Remark} The statement follows from the proof of theorem \ref{th_infinitizemal}.

Let us show a sufficient condition for the function   $(c_1,c_2):\ir\rightarrow \mathbb{R}^2$ to provide a Nash equilibrium. Denote $$H_1(t,x,s)\triangleq \max_{u\in P}\min_{v\in Q}\langle s,f(t,x,u,v)\rangle,\ \ H_2(t,x,s)\triangleq \max_{v\in Q}\min_{u\in P}\langle s,f(t,x,u,v)\rangle. $$ Let $(c_1,c_2):\ir\rightarrow\mathbb{R}^2$, $(t,x)\in \ir$, $w\in\mathbb{R}^n$ define a modulus derivative at the position $(t,x)$ in the direction $w\in\mathbb{R}^n$ by the rule
$${\rm d}_{abs}(c_1,c_2)(t,x;w)\triangleq \liminf_{\delta\downarrow 0,w'\rightarrow w}\frac{|c_1(t+\delta,x+\delta w')-c_1(t,x)|+|c_2(t+\delta,x+\delta w')-c_2(t,x)|}{\delta}. $$

\begin{coll}\label{cor_gen_HJ}
Suppose that the function $(c_1,c_2):\ir\rightarrow\mathbb{R}^2$ is continuous, $(c_1(\vartheta_0,\cdot),c_2(\vartheta_0,\cdot))=(\sigma_1(\cdot),\sigma_2(\cdot))$, for each $i$ the function $c_i$ is upper viscosity solution of the equation \begin{equation}\label{HJ_aux}\frac{\partial c_i}{\partial t}+H_i(t,x,\nabla c_i)=0, \end{equation}
 and for all $(t,x)\in \ir$ $$\inf_{w\in\mathcal{F}(t,x)}{\rm d}_{abs}(c_1,c_2)(t,x;w)=0. $$ Then for all $(t,x)\in \ir$ the pair of numbers  $(c_1(t,x),c_2(t,x))$ is a Nash equilibrium payoff of the game.
\end{coll}

Corollary  \ref{cor_gen_HJ} follows from the definition of modulus derivative and the property of upper solution of equation  (\ref{HJ_aux}) \cite{Subb_book}: $\omega_i(t,x)\leq c_i(t,x)$ for all $(t,x)\in\ir$.

Let us show that the suggested method is a generalization of the method based on the system of Hamilton-Jacobi equations. This method provide a Nash solution in the class of continuous strategies \cite{Basar}.
\begin{Pl}\label{pl_hj}
Let the function $(\varphi_1,\varphi_2):\ir\rightarrow\mathbb{R}^2$ be differentiable, and $(\varphi_1(\vartheta_0,\cdot),\varphi_2(\vartheta_0,\cdot))=(\sigma_1(\cdot),\sigma_2(\cdot))$. Suppose that the function $(\varphi_1,\varphi_2)$ satisfies the following condition: for all positions $(t,x)\in \ir$ there exist $u^n\in P$, $v^n\in Q$ such that
\begin{equation}\label{HJ_u_star_def}
 \max_{u\in P}\big\langle\nabla\varphi_1(t,x),f(t,x,u,v^n)\big\rangle= \big\langle\nabla\varphi_1(t,x),f(t,x,u^n,v^n)\big\rangle, \end{equation}
\begin{equation}\label{HJ_v_star_def}\max_{v\in Q}\big\langle\nabla\varphi_2(t,x),f(t,x,u^n,v)\big\rangle=\big \langle\nabla\varphi_2(t,x),f(t,x,u^n,v^n)\big\rangle \end{equation} и
\begin{equation}\label{HJ_eq_varphi}
\frac{\partial \varphi_i(t,x)}{\partial t}+\big\langle\nabla\varphi_i(t,x),f(t,x,u^n,v^n)\big\rangle=0,\ \ i=1,2.
\end{equation}
Then the function $(\varphi_1,\varphi_2)$ satisfies the conditions of corollary~$\ref{cor_gen_HJ}$.
\end{Pl}

This proposition is proved in section~\ref{sec_weak}.

If one can choose the pair $(u^n,v^n)$ for each position  $(t,x)\in \ir$ and the pair of directions $s_1,s_2\in\mathbb{R}^n$ uniquely, then the Hamiltonians  $\mathcal{H}_i$ are well defined by the rule
$$\mathcal{H}_i(t,x,s_1,s_2)\triangleq \big\langle s_i,f(t,x,u^n,v^n)\big\rangle, \ \ i=1,2. $$ In this case condition (\ref{HJ_eq_varphi}) is equal to the following one: $(\varphi_1,\varphi_2)$ is a solution of the system
$$\frac{\partial \varphi_i}{\partial t}+\mathcal{H}_i(t,x,\nabla\varphi_1,\nabla\varphi_2)=0,\ \ i=1,2. $$

\section{Example}

Consider the nonzero-sum differential game with the dynamic \begin{equation}\label{sys_example}
    \left\{
    \begin{array}{ccc}
      \dot{x} & = & u,\\
      \dot{y} & = & v,
    \end{array}
    \right.
\end{equation} $t\in [0,1]$, $u,v\in [-1,1]$. Payoffs are determined by the formulas $\sigma_1(x,y)\triangleq -|x-y|$, $\sigma_2(x,y)\triangleq y$. We recall that each player wants to maximize his payoff.

In order to determine the multivalued map  $\mathcal{N}:[0,1]\times \mathbb{R}^2\rightarrow\mathcal{P}(\mathbb{R}^2)$, we determine auxiliary multivalued maps $\mathcal{S}_i:[0,1]\times \mathbb{R}^2\rightarrow\mathcal{P}(\mathbb{R})$, such that $$\mathcal{S}_i(t,x_*,y_*) \triangleq\big\{z\in\mathbb{R}: \omega_i(t,x_*,y_*)\leq z\leq c_i^+(t,x_*,y_*)\big\}.$$ Here $$c_i^+(t,x,y)\triangleq \sup_{u\in \mathcal{U},v\in\mathcal{V}}\sigma_i\Bigg(x+\int\limits_{t}^{\vartheta_0}u(\xi)d\xi,y+\int\limits_{t}^{\vartheta_0}v(\xi)d\xi\Bigg). $$  Obviously,
\begin{equation}\label{N_S_subsets}
\mathcal{N}(t,x_*,y_*)\subset \mathcal{S}_1(t,x_*,y_*)\times\mathcal{S}_2(t,x_*,y_*).
\end{equation}

First we determine the map $\mathcal{S}_2$. The value function of the game $\Gamma_2$  is equal to $\omega_2(t,x_*,y_*)=y_*+(1-t)$. Also, $c_2^+(t,x_*,y_*)=y_*+(1-t)$. Consequently
\begin{equation}\label{S_2_formula}
\mathcal{S}_2(t,x_*,y_*)=y_*+(1-t).
\end{equation}

Let us determine the set $\mathcal{S}_1$.  Programmed iteration method \cite{Math_sb} yields that   $$\omega_1(t,x_*,y_*)=-|x_*-y_*|. $$ Moreover $$c^+_1(t,x_*,y_*)=\min\big\{-|x_*-y_*|+2(1-t),0\big\}. $$ We obtain that
\begin{equation}\label{S_1_formula}
\mathcal{S}_1(t,x_*,y_*)=\big[\omega_1(t,x_*,y_*),c^+_1(t,x_*,y_*)\big].
\end{equation}

Now we compute the map $\mathcal{N}(t,x_*,y_*)$. The linearity of right hand of  (\ref{sys_example}) and the convexity of restrictions on control yield that any control can be substitute by the pair of constant controls $(u,v)\in P\times Q$. We have that for all $(J_1,J_2)\in \mathcal{S}_1(t,x_*,y_*)\times S_2(t,x_*,y_*)$ $$\Dh{\mathcal{N}}{t,x_*,y_*;(J_1,J_2),(u,v)}\\\geq \liminf_{\delta\downarrow 0, v'\rightarrow v}\frac{|y_*+\delta v'+(1-t-\delta)-y_*-(1-t)|}{\delta}=|v-1|. $$ Therefore, if  $\Dh{\mathcal{N}}{t,x_*,y_*;(J_1,J_2),w}=0$ for the pair $w=(u,v)$, then $v=1$.

First we consider the case $y_*\geq x_*$. Let $(J_1,J_2)\in\mathcal{N}(t,x_*,y_*)$. There exists a motion $(x(\cdot),y(\cdot))\in {\rm Sol}(t,x_*,y_*)$ such that $(J_1,J_2)\in \mathcal{N}(\theta,x(\theta),y(\theta))$, $\theta\in [t,1]$. Since $\Dh{\mathcal{N}}{t,x_*,y_*;(J_1,J_2),w}=0$ only if $v=1$,  there exists $u\in [-1,1]$ such that $x(1)=x_*+u(1-t)$, $y(1)=y_*+(1-t)$. Consequently $y(1)\geq x(1)$. From condition (N2) we get that
\begin{multline*}J_1=\sigma_1(x(1),y(1))=-y(1)+x(1)=-y_*-(1-t)+x_*+u(1-t)\\ =-y_*+x_*+(1-t)(u-1)\leq -y_*+x_*=-|x_*-y_*|. \end{multline*} The equality is achieved only if $u=1$. Condition (N1) yields that the the following inclusion is fulfilled $$\mathcal{N}(t,x_*,y_*)\subset \{(-|x_*-y_*|,y_*+(1-t))\}. $$ Substituting value $(1,1)$ for $w$ in formula for $\Dh{\mathcal{N}}{t,x_*,y_*;(-|x_*-y_*|,y_*+(1-t)),w}$ we claim that for $y_*\geq x_*$
$$
 \mathcal{N}(t,x_*,y_*)= \big\{(-|x_*-y_*|,y_*+(1-t))\big\}.
$$

Now let $y_*<x_*$. We shall show that
\begin{multline*}
\mathcal{N}(t,x_*,y_*)= \mathcal{S}_1(t,x_*,y_*)\times\mathcal{S}_2(t,x_*,y_*)=\\\Big[-|x_*-y_*|,\min\big\{-|x_*-y_*|+2(1-t),0\big\}\Big]\times \big\{y_*+(1-t)\big\}.
\end{multline*}
Clearly, conditions (N1) and (N2) hold for this map. Let $\gamma_0$ be a maximal number of segment  $[0,2]$ such that $-|x_*-y_*|+\gamma_0(1-t)\leq 0$. If $(J_1,J_2)\in \mathcal{N}(t,x_*,y_*)$, then $J_2=y_*+(1-t)$, $J_1=-|x_*-y_*|+d(1-t)$ for some $d\in [0,\gamma_0]$. Let us prove that there exists a number $\delta>0$ with the property \begin{equation}\label{j_1_inclus}
(J_1,J_2)\in\mathcal{N}(t+\delta,x_*+\delta u,y_*+\delta)
\end{equation}
for $u=1-d$. It is sufficient to prove that $$J_1\in \Big[y_*-x_*+\delta d,\min\big\{y_*-x_*+\delta d+2(1-t-\delta),0\big\}\Big].$$ Indeed, $y_*-x_*+d(1-t)\geq y_*-x_*+\delta d$ for $\delta<(1-t)$.
Since $d\leq\gamma_0$, we obtain that
$$J_1=y_*-x_*+d(1-t)\leq y_*-x_*+\delta d+\gamma_0(1-t-\delta). $$
Also $$y_*-x_*+\delta d+\gamma_0(1-t-\delta)\leq\min\big\{y_*-x_*+d\delta+2(1-t-\delta);0\big\}.$$ Actually, since $\gamma_0\leq 2$, the following inequality is fulfilled $$y_*-x_*+\delta d+\gamma_0(1-t-\delta)\leq y_*-x_*+d\delta+2(1-t-\delta).$$ Moreover  $y_*-x_*+\delta d+\gamma_0(1-t-\delta)\leq\delta d-\gamma_0 \delta\leq 0$.
 Thus the condition
$$J_1\leq  \min\big\{y_*-x_*+\delta d+2(1-t-\delta),0\big\}=y_*-x_*+\delta d+\gamma_1(1-t-\delta)$$ is valid also.
It follows from (\ref{j_1_inclus}) that $$\Dh{\mathcal{N}}{t,x_*,y_*;(J_1,J_2),(1-d,1)}=0. $$ Since $\mathcal{N}(t,x_*,y_*)$ coincide with the set $\mathcal{S}_1(t,x_*,y_*)\times \mathcal{S}_2(t,x_*,y_*)$ in this case, we claim that the set $\mathcal{N}(t,x_*,y_*)$ is  Nash value of the game at the position $(t,x_*,y_*)$.

Let us compare the obtained result with the method based on system of Hamilton-Jacobi equations \cite{Basar}. In considered case the system of equations is given by
\begin{equation}\label{hj_sys}
\left\{
\begin{array}{ll}
\displaystyle\frac{\partial \varphi_1}{\partial t}+\displaystyle\frac{\partial \varphi_1}{\partial x} \hat{u}(t,x,y) +\displaystyle\frac{\partial \varphi_1}{\partial y} \hat{v}(t,x,y) & =0 \\[3ex]
\displaystyle\frac{\partial \varphi_2}{\partial t}+\displaystyle\frac{\partial \varphi_2}{\partial x} \hat{u}(t,x,y) +\displaystyle\frac{\partial \varphi_2}{\partial y} \hat{v}(t,x,y) & =0.
\end{array}
\right.
\end{equation} Here the values $\hat{u}(t,x,y)$ and $\hat{v}(t,x,y)$ are determined by the following conditions
$$\frac{\partial \varphi_1(t,x,y)}{\partial x} \hat{u}(t,x,y)=\max_{u\in P}\Big[\frac{\partial \varphi_1(t,x,y)}{\partial x}u\bigg], \ \ \frac{\partial \varphi_2(t,x,y)}{\partial y} \hat{v}(t,x,y)=\max_{v\in Q}\bigg[\frac{\partial \varphi_2(t,x,y)}{\partial y}v\bigg].$$
It follows from Proposition \ref{pl_hj} that if a pair of functions  $(\varphi_1,\varphi_2)$ is a solution of the system (\ref{hj_sys}), then $\varphi_2(t,x,y)=y+(1-t)$. Thus $v_*(t,x,y)=1$. Consequently, the system~(\ref{hj_sys}) reduces to the equation
\begin{equation}\label{hj_sys_one}
 \frac{\partial \varphi_1}{\partial t}+\left|\frac{\partial \varphi_1}{\partial x}\right|+\frac{\partial \varphi_1}{\partial y}=0.
\end{equation} By \cite[Theorem 5.6]{Subb_book} we obtain that the function $$
\varphi_1(t,x,y)=\left\{
\begin{array}{ll}
x-y, & x\leq y, \\
-x+y+2(1-t), & x>y, -x+y+2(1-t)<0,\\
0, & x>y, -x+y+2(1-t)\geq 0
 \end{array}
\right.
 $$ is a minimax solution of equation (\ref{hj_sys_one}). Indeed if $\varphi_1$ is smooth at $(t,x,y)$ then equation (\ref{hj_sys_one}) is fulfilled in classical sense. On a planes $\{(t,x,y):x=y\}$, $\{(t,x,y):-x+y+2(1-t)=0\}$ we have that the Clarke subdifferential is the convex hull of two limit of partial derivatives of  the function $\varphi_1$. By well-known properties of subdifferentials and superdifferentail, the continuity  and positive homogeneity of equation (\ref{hj_sys_one}) we obtain that $\varphi_1$ satisfies conditions U4 and L4 of \cite{Subb_book}.

The function $\varphi_1$ is nonsmooth. Since the minimax solution is unique, and any classical solution is minimax, we claim that system (\ref{hj_sys}) have no classical solution.
 One may obtain from the formulae for $\mathcal{N}(t,x,y)$ that $(\varphi_1(t,x,y),\varphi_2(t,x,y))\in\mathcal{N}(t,x,y)$. Moreover, $$\varphi_1(t,x,y)=\max\Big\{J_1\in\mathbb{R}:\exists J_2\in\mathbb{R}\ \ (J_1,J_2)\in \mathcal{N}(t,x,y)\Big\},$$ $$\{\varphi_2(t,x,y)\}=\Big\{J_1\in\mathbb{R}:\exists J_2 \ \ (J_1,J_2)\in\mathcal{N}(t,x,y)\Big\}. $$

In other words, the value $(\varphi_1(t,x,y),\varphi_2(t,x,y))$ is the maximal Nash equilibrium payoff of the game at the position $(t,x,y)$.

One can check that the pair of functions $(\varphi_1,\varphi_2)$ satisfies the conditions of corollary \ref{cor_gen_HJ}. Simultaneity,  there exists a family of function satisfying the condition of corollary \ref{cor_gen_HJ}. Actually, if $\gamma\in [0,2]$, then put
  $$c_1^\gamma(t,x_*,y_*)=\left\{
\begin{array}{ll}
 -|x_*-y_*|, & y_*\geq x_*; \\[1ex]
 \min\big\{-|x_*-y_*|+\gamma(1-t);0\big\}, & y_*<x_*
\end{array}
\right.$$
$$ c^\gamma_2(t,x_*,y_*)=y_*+(1-t).$$
Let us show that the pair of functions $(c_1^\gamma,c_2^\gamma)$ satisfy the conditions of corollary \ref{cor_gen_HJ}.
We have that in our case $$H_1(t,x,y,s_x,s_y)=|s_x|-|s_y|, \ \ H_2(t,x,y,s_x,s_y)=|s_y|-|s_x|. $$
First we prove that the functions $c_i^\gamma$ are the upper solution of equations (\ref{HJ_aux}). By  \cite[condition U4]{Subb_book} it  suffices to show that for all $(t,x,y)\in [t_0,\vartheta_0]\times\mathbb{R}^2$ $(a,s_x,s_y)\in D^-c_i^\gamma(t,x,y)$ the following inequalities holds
\begin{equation}\label{u_4_example}
a+H_i(s_x,s_y)\leq 0,\ \  i=1,2.
\end{equation}
Here $D^-$ denotes the subdifferential  \cite[(6.10)]{Subb_book}.
The computing of subdifferentials   gives that
$$D^-c_1^\gamma(t,x,y)=\left\{\begin{array}{ll}
  \{(0,-1,1)\}, & y>x, \\[0.5ex]
  \{(0,0,0)\}, & y<x<y+\gamma(1-t), \\[0.5ex]
  \{(-\gamma,-1,1)\}, & x>y+\gamma(1-t), \\[0.5ex]
  \{(0,\lambda,-\lambda):\lambda\in [0,1]\}, & x=y, \\[0.5ex]
  \{(-\lambda\gamma,-\lambda,\lambda):\lambda\in [0,1]\}, & x=y+\gamma(1-t).
\end{array}\right.$$
$$D^-c_2^\gamma(t,x,y)=\{(-1,0,1)\}. $$ Substituting the values of subdifferentials, we get that (\ref{u_4_example}) is valid for $i=1,2$.

Also $c_i(1,x_*,y_*)=\sigma_i(x_*,y_*)$.
Moreover, ${\rm d}_{abs}(c_1^\gamma,c_2^\gamma)(t,x_*,y_*;1-d,1)=0$ for
$$d=\left\{
\begin{array}{ll}
 0, & y_*\geq x_*; \\[0.5ex]
 \max\Big\{r\in [0,\gamma]:-|x_*-y_*|+r(1-t)\leq 0\Big\}, & y_*<x_*.
\end{array}
\right. $$

Note that $(\varphi_1,\varphi_2)=(c_1^2,c_2^2)$.

\section{Weak invariance of the set of values}\label{sec_weak}

In this section the statements formulated in section \ref{sec_main} are proved.

\proof[Proof of Theorem \ref{th_weak_inv}]{
If condition (N3) holds, then one can put $\theta=\vartheta_0$.

Now suppose that for all $(t_*,x_*)$, $(J_1,J_2)\in \mathcal{T}(t_*,x_*)$ there exist $\theta\in [t_*,\vartheta_0]$ and a motion
$y(\cdot)\in {\rm Sol}(t_*,x_*)$, such that the following condition is fulfilled \begin{equation}\label{theta_def}(J_1,J_2)\in\mathcal{T}(t,y(t)),\ \  t\in [t_*,\theta].\end{equation}

Let $\Theta$ be a set of moments $\theta$ satisfying condition (\ref{theta_def}) for some $y(\cdot)\in {\rm Sol}(t_*,x_*)$. Denote $\tau\triangleq \sup\Theta$. We have that $\tau\in \Theta$. Indeed, let a sequence $\{\theta_k\}_{k=1}^\infty\subset\Theta$ tend to $\tau$. One can assume that $\theta_k<\theta_{k+1}\leq\tau$. For every $k$ condition (\ref{theta_def}) is valid under $\theta=\theta_k$, $y(\cdot)=y_k(\cdot)\in {\rm Sol}(t_*,x_*)$. The compactness of bundle of motions yields that $y_k(\cdot)\rightarrow y^*(\cdot)$, $k\rightarrow\infty$, here $y^*(\cdot)$ is an element of ${\rm Sol}(t_*,x_*)$. The closeness of the map $\mathcal{T}$ gives that $(J_1,J_2)\in \mathcal{T}(\theta_k,y^*(\theta_k))$. By the same argument we claim that $(J_1,J_2)\in \mathcal{T}(\tau,y^*(\tau))$. Denote $x^*=y^*(\tau)$.

Let us show that $\tau=\vartheta_0$. If $\tau<\vartheta_0$, then there exist a motion $\hat{y}(\cdot)\in{\rm Sol}(\tau,x^*)$ and a moment  $\theta'>\tau$ such that $(J_1,J_2)\in \mathcal{T}(t,\hat{y}(t))$, $t\in [\tau,\theta']$.
Consider a motion
$$\tilde{y}(t)\triangleq \left\{\begin{array}{cc}
 y^*(t), & t\in [t_*,\tau],\\[0.5ex]
 \hat{y}(t), & t\in [\tau,\theta'].
\end{array}\right.
 $$  By definition of $\theta'$ it follows that (\ref{theta_def}) is valid under $\theta=\theta'$, $y(\cdot)=\tilde{y}(\cdot)$. Thus $\theta'\in\Theta$, this contradicts with the choice of $\tau$. Consequently, $\tau=\vartheta_0$ and condition (N3) holds.
}

\proof[Proof of theorem \ref{th_infinitizemal}]{ Let us introduce a graph of map $\mathcal{T}$ $${\rm gr}\mathcal{T}\triangleq \{(t,x,J_1,J_2):(t,x)\in \ir, (J_1,J_2)\in\mathcal{T}(t,x)\}.$$

One can reformulate the condition of theorem \ref{th_weak_inv} in the following way: the graph of $\mathcal{T}$ is weakly invariant under the differential inclusion $$\left(\begin{array}{c}\dot{x}\\
\dot{J}_1\\
\dot{J}_2
\end{array}\right)
\in \widehat{\mathcal{F}}(t,x)\triangleq {\rm co}\left\{\left(\begin{array}{c}
 f(t,x,u,v)\\
 0\\
 0
\end{array}\right):u\in P, v\in Q
\right\}. $$

The condition of weak invariance of the multivalued map $\mathcal{T}$ under differential inclusion~$\widehat{\mathcal{F}}$ is equivalent \cite{Subb_book, Subbotin_invar} to the following condition
\begin{equation}\label{Dt_nonempty}D_t({\rm gr}\mathcal{T})(t,x,J_1,J_2)\cap \widehat{\mathcal{F}}(t,x)\neq \varnothing \end{equation}
for all $(t,x)\in \ir$, $(J_1,J_2)\in\mathcal{T}(t,x)$. Here $D_t$ denotes right-hand derivative in $t$. Let $\mathcal{G}\subset [t_0,\vartheta_0]\times \mathbb{R}^m$, $\mathcal{G}[t]$ denote a section of $G$ by $t$: $$\mathcal{G}[t]\triangleq \{w\in \mathbb{R}^m:(t,x)\in \mathcal{G}\},$$ symbol ${\rm d}$ denote Euclidian distance between a point and a set. Following \cite{Subb_book}, \cite{Subbotin_invar} put $$(D_t\mathcal{G})(t,y)\triangleq \left\{h\in\mathbb{R}^m: \liminf_{\delta\rightarrow 0}\frac{{\rm d}(y+\delta h; \mathcal{G}[t+\delta])}{\delta}=0\right\}. $$

Let us show that conditions (\ref{dh_zero}) and (\ref{Dt_nonempty}) are equivalent.

Condition (\ref{dh_zero}) means that for every pair $(J_1,J_2)\in\mathcal{T}(t,x)$ the following condition holds: $$\inf_{w\in\mathcal{F}(t,x)}\liminf_{\delta\downarrow 0,\gamma\in\mathbb{R}^n,\|\gamma\|\downarrow 0} \frac{{\rm dist}\big[(J_1,J_2),\mathcal{T}(t+\delta,x+\delta(w+\gamma))\big]}{\delta}=0. $$

The lower boundary by $w$ in the formula
$$
\inf_{w\in\mathcal{F}(t,x)}\liminf_{\delta\downarrow 0,\gamma\in\mathbb{R}^n,\|\gamma\|\downarrow 0}\frac{{\rm dist}\big[(J_1,J_2),\mathcal{T}(t+\delta,x+\delta(w+\gamma))\big]}{\delta}$$ is attained for all $(J_1,J_2)\in\mathcal{T}(t,x)$. Indeed, let $\{w^r\}_{r=1}^\infty$ be a minimizing sequence. By the compactness of $\mathcal{F}(t,x)$ one can assume that $w^r\rightarrow w^*$, $r\rightarrow \infty$, $w^*\in\mathcal{F}(t,x)$. Let us show that
\begin{multline}\label{inf_realize}
\tilde{b}\triangleq \inf_{w\in\mathcal{F}(t,x)}\liminf_{\delta\downarrow 0,\gamma\in\mathbb{R}^n,\|\gamma\|\downarrow 0} \frac{{\rm dist}\big[(J_1,J_2),\mathcal{T}(t+\delta,x+\delta(w+\gamma))\big]}{\delta}\\=
\liminf_{\delta\downarrow 0,\gamma\in\mathbb{R}^n,\|\gamma\|\downarrow 0} \frac{{\rm dist}\big[(J_1,J_2),\mathcal{T}(t+\delta,x+\delta(w^*+\gamma))\big]}{\delta}.
\end{multline} Indeed for every $r\in\mathbb{N}$ there exist a sequences $\{\delta^{r,k}\}_{k=1}^\infty$, $\{\gamma^{r,k}\}_{k=1}^\infty$ such that $\delta^{r,k}, \|\gamma^{r,k}\|\rightarrow 0$, $k\rightarrow \infty$ and
\begin{multline*}
b^r\triangleq \liminf_{\delta\downarrow 0,\gamma\in\mathbb{R}^n,\|\gamma\|\downarrow 0} \frac{{\rm dist}\big[(J_1,J_2),\mathcal{T}(t+\delta,t+\delta(w^r+\gamma))\big]}{\delta}\\=
\lim_{k\rightarrow \infty}\frac{{\rm dist}\big[(J_1,J_2),\mathcal{T}(t+\delta^{r,k},t+\delta^{r,k}(w^r+\gamma^{r,k}))\big]}{\delta^{r,k}}.
\end{multline*}

Let $\hat{k}(r)$ be a number such that
$$\delta^{r,\hat{k}(r)}, \|\gamma^{r,\hat{k}(r)}\|, \left|\frac{{\rm dist}\big[(J_1,J_2),\mathcal{T}(t+\delta^{r,\hat{k}(r)},t+\delta^{r,\hat{k}(r)}(w^r+\gamma^{r,\hat{k}(r)}))\big]}{\delta^{r,\hat{k}(r)}}-b^r\right|<2^{-r}.$$
Put $\hat{\delta}^r\triangleq \delta^{r,\hat{k}(r)}$, $\hat{\gamma}^r\triangleq \gamma^{r,\hat{k}(r)}+w^r-w^*$. Note that $\hat{\delta}^r, \|\hat{\gamma}^r\|\rightarrow 0$, $r\rightarrow \infty$.

We have that
\begin{multline}\label{inf_ineq_r}
 \inf_{w\in\mathcal{F}(t,x)}\liminf_{\delta\downarrow 0,\gamma\in\mathbb{R}^n,\|\gamma\|\downarrow 0} \frac{{\rm dist}\big[(J_1,J_2),\mathcal{T}(t+\delta,x+\delta(w+\gamma))\big]}{\delta} \\\leq
\liminf_{\delta\downarrow 0,\gamma\in\mathbb{R}^n,\|\gamma\|\downarrow 0} \frac{{\rm dist}\big[(J_1,J_2),\mathcal{T}(t+\delta,x+\delta(w^*+\gamma))\big]}{\delta}\\\leq \lim_{r\rightarrow \infty}\frac{{\rm dist}\big[(J_1,J_2),\mathcal{T}(t+\hat{\delta}^r,x+\hat{\delta}^r(w^*+\hat{\gamma}^r))\big]}{\hat{\delta}^r}.
\end{multline}

Further, \begin{multline*}
\frac{{\rm dist}\big[(J_1,J_2),\mathcal{T}(t+\hat{\delta}^r,x+\hat{\delta}^r(w^*+\hat{\gamma}^r))\big]}{\hat{\delta}^r}\\=
\frac{{\rm dist}\big[(J_1,J_2),\mathcal{T}(t+\delta^{r,\hat{k}(r)},x+\delta^{r,\hat{k}(r)}(w^*+\gamma^{r,\hat{k}(r)}+w^r-w^*))\big]}{\delta^{r,\hat{k}(r)}}\\=\frac{{\rm dist}\big[(J_1,J_2),\mathcal{T}(t+\delta^{r,\hat{k}(r)},x+\delta^{r,\hat{k}(r)}(w^k+\gamma^{r,\hat{k}(r)}))\big]}{\delta^{r,\hat{k}(r)}}\\\leq b^r+2^{-r}\rightarrow \tilde{b}, \ \ r\rightarrow \infty.
       \end{multline*}
We have that in (\ref{inf_ineq_r}) right and left hands are equal. This means that condition (\ref{inf_realize}) is valid.

Thus, condition (\ref{dh_zero}) is equivalent to the following one: for all $(J_1,J_2)\in\mathcal{T}(t,x)$ there exists $w\in\mathcal{F}(t,x)$ such that
\begin{multline}\label{diff_zero}
\liminf_{\delta\downarrow 0,\gamma\in\mathbb{R}^n,\|\gamma\|\downarrow 0}\frac{{\rm dist}\big[(J_1,J_2),\mathcal{T}(t+\delta,x+\delta(w+\gamma))\big]}{\delta}=\\
 \liminf_{\delta\downarrow 0,\gamma\in\mathbb{R}^n,\|\gamma\|\downarrow 0} \inf\left\{\frac{|\zeta_1-J_1|+|\zeta_2-J_2|}{\delta}:(\zeta_1,\zeta_2)\in \mathcal{T}(t+\delta,x+\delta (w+\gamma)) \right\}=0.
\end{multline}

Now let us prove that this condition is equivalent to condition (\ref{Dt_nonempty}).

First we assume that condition (\ref{Dt_nonempty}) is valid. This means that there exist sequences $\{\delta^k\}_{k=1}^\infty\subset\mathbb{R},$ $\{\gamma^k\}_{k=1}^\infty\subset\mathbb{R}^n,$ $\{\varepsilon_1^k\}_{k=1}^\infty, \{\varepsilon_2^k\}_{k=1}^\infty\subset\mathbb{R}$ such that
\begin{itemize}
\item $\delta^k,\|\gamma^k\|, \varepsilon^k_1, \varepsilon^k_2\rightarrow 0,$ $k\rightarrow \infty$;
\item $\Big(t+\delta^k,x+\delta^k(w+\gamma^k), J_1+\delta^k\varepsilon_1^k, J_2+\delta^k\varepsilon_2^k\Big)\in {\rm gr}\mathcal{T}$.
\end{itemize}

One can reformulate the second condition as $$(J_1+\delta^k\varepsilon_1^k,J_2+\delta^k\varepsilon_2^k)\in \mathcal{T}(t+\delta^k,t+\delta^k(w+\gamma^k)). $$
Thus, $$\inf\left\{\frac{|\zeta_1-J_1|+|\zeta_2-J_2|}{\delta^k}:(\zeta_1,\zeta_2)\in \mathcal{T}(t+\delta^k,x+\delta (w+\gamma^k)) \right\}=\varepsilon^k_1+\varepsilon^k_2. $$ By the choice $\{\varepsilon_1^k\}$, $\{\varepsilon_2^k\}$ we obtain that condition (\ref{diff_zero}) holds.

Now let condition (\ref{diff_zero}) be fulfilled, prove that (\ref{Dt_nonempty}) is valid. Indeed, let $\{\delta^k\}_{k=1}^\infty$, $\{\gamma\}_{k=1}^\infty$ be a minimizing sequence. By compactness of the sets  $\mathcal{T}\big(t+\delta^k,x+\delta^k(w+\gamma^k)\big)$ for each $k$ there exist $\varepsilon^k_1$ and $\varepsilon^k_2$ such that
$$\big(J_1+\delta^k\varepsilon^k_1,J_2+\delta^k\varepsilon^k_2\big)\in\mathcal{T}(t,x+\delta^k(w+\gamma^k)). $$

It follows from (\ref{diff_zero}) that $\varepsilon^k_1,\varepsilon^k_2\rightarrow 0$, $k\rightarrow \infty$. Let us estimate
${{\rm d}\Big((w+\delta^kw,J_1,J_2),{\rm gr}\mathcal{T}[t+\delta^k]\Big).}$ We have that $$\Big(t+\delta^k,x+\delta^k(w+\gamma^k), J_1+\delta^k\varepsilon^k_1,J_2+\delta^k\varepsilon^k_2\Big)\in {\rm gr}\mathcal{T}.$$ Consequently,
$${\rm d}\Big((w+\delta^kw,J_1,J_2),{\rm gr}\mathcal{T}[t+\delta^k]\Big)\leq \delta^k\sqrt{\|\gamma^k\|^2+(\varepsilon^k_1)^2+(\varepsilon^k_2)^2}. $$ The convergence $\delta^k,\|\gamma^k\|, \varepsilon^k_1,\varepsilon^k_2\rightarrow 0$ as $k\rightarrow \infty$ yields the equality $$\left(\begin{array}{c}w\\
0\\
0\end{array}\right)\in D_t({\rm gr}\mathcal{T})(t,x,J_1,J_2). $$ Since $\widehat{\mathcal{F}}(t,x)=\mathcal{F}\times \{(0,0)\}$, we claim that (\ref{Dt_nonempty}) is fulfilled.

}

\proof[Proof of Proposition \ref{pl_hj}]{
It follows from (\ref{HJ_u_star_def}) and the Isaacs condition that
$$\big\langle\nabla\varphi_1(t,x),f(t,x,u^n,v^n)\big\rangle \geq \max_{u\in P}\min_{v\in Q}\big\langle \nabla\varphi_1(t,x),f(t,x,u,v)\big\rangle=H_1(t,x,\nabla\varphi_1(t,x)).$$ Analogously, it follows from (\ref{HJ_v_star_def}) and the Isaacs condition that $$\big\langle\nabla\varphi_2(t,x),f(t,x,u^n,v^n)\big\rangle \geq \max_{v\in Q}\min_{u\in P}\big\langle \nabla\varphi_2(t,x),f(t,x,u,v)\big\rangle=H_2(t,x,\nabla\varphi_2(t,x)).$$ Therefore, using (\ref{HJ_eq_varphi}) we claim that $$\frac{\partial\varphi_i(t,x)}{\partial t}+H_i(t,x,\nabla\varphi_i(t,x))\leq 0, \ \ i=1,2. $$ Since the function $\varphi_i$ is differentiable, its subdifferential at the position $(t,x)$ is equal to $\{\partial\varphi_1(t,x)/\partial t,\nabla\varphi_1(t,x)\}$. Consequently, the function $\varphi_1$ is the upper solution of equation (\ref{HJ_aux}) for $i=1$ \cite[Condition (U4)]{Subb_book}. Analogously, the function $\varphi_2$ is the upper solution of equation (\ref{HJ_aux}) for $i=2$.

Now let us show that ${\rm d}_{abs}(\varphi_1,\varphi_2)(t,x;w)=0$ for $w\in\mathcal{F}(t,x)$. Put $w=f(t,x,u^n,v^n)$. Indeed,
\begin{multline*}
{\rm d}_{abs}(\varphi_1,\varphi_2)(t,x;w)\\=\liminf_{\delta\downarrow 0,\|\gamma\|\rightarrow 0}\frac{|\varphi_1(t+\delta,x+\delta(w+\gamma))-\varphi_1(t,x)|+|\varphi_2(t+\delta,x+\delta(w+\gamma))-\varphi_2(t,x)|}{\delta}.
\end{multline*} Let $\{\delta^k\}_{k=1}^\infty\subset\mathbb{R}$, $\{\gamma^k\}_{k=1}^\infty\subset\mathbb{R}^n$ be a minimizing sequence. Then
\begin{multline*}
{\rm d}_{abs}(\varphi_1,\varphi_2)(t,x;w)\\=\lim_{k\rightarrow \infty}\frac{|\varphi_1(t+\delta^k,x+\delta^k(w+\gamma^k))-\varphi_1(t,x)|+|\varphi_2(t+\delta^k,x+
\delta(w+\gamma^k))-\varphi_2(t,x)|}{\delta^k}\\=
\lim_{k\rightarrow \infty}\frac{1}{\delta^k}\Bigl[\Bigl|\frac{\partial\varphi_1(t,x)}{\partial t}\delta^k+\langle\nabla\varphi_1(t,x),\delta^k(w+\gamma^k)\rangle+o(\delta^k) \Bigr|\\+\Bigl|\frac{\partial\varphi_2(t,x)}{\partial t}\delta^k+\langle\nabla\varphi_2(t,x),\delta^k(w+\gamma^k)\rangle+o(\delta^k) \Bigr|\Bigr]\\=\Bigl|\frac{\partial\varphi_1(t,x)}{\partial t}+\langle\nabla\varphi_1(t,x),w\rangle\Bigr|+\Bigl|\frac{\partial\varphi_2(t,x)}{\partial t}+\langle\nabla\varphi_2(t,x),w\rangle\Bigr|.
\end{multline*} By choice of $w=f(t,x,u^n,v^n)$ and condition (\ref{HJ_eq_varphi}) we have that $$\frac{\partial\varphi_1(t,x)}{\partial t}+\langle\nabla\varphi_1(t,x),w\rangle=\frac{\partial\varphi_2(t,x)}{\partial t}+\langle\nabla\varphi_2(t,x),w\rangle=0. $$ Thus ${\rm d}_{abs}(\varphi_1,\varphi_2)(t,x;w)=0$.

}

\end{document}